\documentclass[11pt]{article}
\usepackage{amsmath}
\usepackage{amssymb}
\def\Z{\mathbb{Z}}
\def\N{\mathbb{N}}
\def\Be{\mathcal{B}}

\newtheorem{theorem}{\hspace*{\parindent}Theorem}
\newtheorem{lemma}{\hspace*{\parindent}Lemma}

\title{A note on Fox's $H$ function in the light of Braaksma's results}
\author{D.B.\:Karp$^{\rm a,b}$\footnote{E-mail: D.B.\:Karp -- \emph{dimkrp@gmail.com}}
\\[10pt]
\small{\textit{$\phantom{1}^a$Far Eastern Federal University, Vladivostok, Russia}}
\\
\small{\textit{$\phantom{1}^b$Institute of Applied Mathematics, FEBRAS, Vladivostok, Russia}}}
\date{}
\begin{document}
\maketitle

\begin{abstract}
In our previous works we found a power series expansion of a particular case of Fox's $H$ function $H^{q,0}_{p,q}$ in a neighborhood of its positive singularity. An inverse factorial series expansion of the integrand of $H^{q,0}_{p,q}$ served as our main tool.  However, a necessary  restriction on parameters is missing in those works. In this note we fill this gap and give a simpler and shorter proof of the expansion around the positive singular point.  We further identify more precisely the abscissa of convergence of the underlying inverse factorial series.  Our  new proof hinges on a slight generalization of a particular case of  Braaksma's theorem about analytic continuation of Fox's $H$ function.
\end{abstract}

\bigskip


\bigskip


\bigskip

\section{Introduction}
Fox's $H$ function is a very general function defined by the Mellin-Barnes integral \eqref{eq:Hmnpq} below.  It was introduced by Fox \cite{Fox} in the context of symmetrical Fourier kernels and has been studied by a number of authors thereafter \cite{Braaksma,KilSaig,MSH,ParKam}.  Fox's $H$ function  has a number of important applications, most notably in statistics
\cite{CoelhoArnold,GuptaTang,MSH} and fractional calculus \cite{KSTBook,MSH}.  Braaksma's mammoth manuscript \cite{Braaksma} remains among the deepest investigations on the topic. In this note we will be only interested in the case when the parameter $\mu$ defined in \eqref{eq:mu-beta-defined} equals zero (Braaksma's work contains a careful study of both cases  $\mu>0$ and $\mu=0$, while the case $\mu<0$ reduces to $\mu>0$ by a simple change of variable).  When $\mu=0$ the integral \eqref{eq:Hmnpq} defining the $H$ function only converges if $|z|<\beta^{-1}$, where $\beta$ is given in \eqref{eq:mu-beta-defined}. It also converges for $|z|>\beta^{-1}$ but over a different integration contour so that the two functions obtained in this way are not, generally speaking, analytic  continuations of each other.  Braaksma constructed analytic continuations of the $H$ function to an infinite set of sectors forming a partition of the Riemann surface of the logarithm.  The points of intersection of the circle $|z|=\beta^{-1}$ with the rays bounding these sectors are singular points of $H$.  He further derived analytic continuation from each such sector to the  domain $|z|>\beta^{-1}$ on this Riemann surface.

Under certain restrictions to be specified below the positive semi-axis $\arg(z)=0$ forms the boundary between two such sectors, so that the point $z=\beta^{-1}$ is a singular point of Fox's $H$ function.  We studied the behavior  of a particular case $H^{q,0}_{p,q}$ of Fox's $H$ function in the neighborhood of this point in our recent papers \cite{KPJMS2018,KPCMFT2017} (as Braaksma's parameter $\mu$ is called $\Delta$ in our papers following \cite{KilSaig}, the case $\mu=0$ was given the denomination ''delta-neutral''). In particular, we found a convergent inverse factorial series for the integrand of the delta-neutral $H^{q,0}_{p,q}$ function and utilized it to derive a power series expansion of this function in a neighborhood of the singular point $z=\beta^{-1}$.  We also deduced recursive and determinantal formulas for the coefficients of this expansion. However, our work contained an error caused by an incorrect citing and use of N{\o}rlund's results in \cite[Theorem~2]{KPJMS2018}.  The consequence of this error was omission of the condition that the scaling factors $\alpha_i$, $\beta_j$ in \eqref{eq:Hmnpq} should all be strictly greater than $1/6$  (and not only positive as claimed) in \cite[Theorem~5]{KPJMS2018} and \cite[Theorem~1]{KPCMFT2017}.

The purpose of this note is to correct, simplify, refine and generalize our results presented in \cite{KPJMS2018,KPCMFT2017}.
We give a simpler proof, refine the convergence abscissa of the inverse factorial series
and present another expansion that works for arbitrary positive scaling factors $\alpha_i$, $\beta_j$. This is done in section~3. In order to obtain these results we needed a slight generalization of a particular case of \cite[Theorem~2]{Braaksma}, which we establish in section~2.  As some parts of Braaksma's work do not seem to be well-understood we also rewrite some of his ideas here.

\section{Braaksma revisited}
As we will need to modify some parts of Braaksma's proof, for reader's convenience we adhere to the notation of \cite{Braaksma}.
Define
\begin{equation}\label{eq:Hmnpq}
H^{m,n}_{p,q}\left(\!z\left|\begin{array}{l}(\alpha_1,a_1),\ldots(\alpha_p,a_p)\\(\beta_1,b_1),\ldots(\beta_q,b_q)\end{array}\right.\!\right)
=H(z)=\frac{1}{2\pi{i}}\int\limits_{C}h(s)z^sds,
\end{equation}
where
$$
z^s=\exp\{s(\log|z|+i\arg{z})\},
$$
and $\arg{z}$ may take any real value. Here
\begin{multline*}
h(s)=\frac{\prod_{j=1}^{n}\Gamma(1-a_j+\alpha_js)\prod_{j=1}^{m}\Gamma(b_j-\beta_js)}{\prod_{j=m+1}^{q}\Gamma(1-b_j+\beta_js)\prod_{j=n+1}^{p}\Gamma(a_j-\alpha_js)}
\\
=\underbrace{\frac{\prod_{j=1}^{p}\Gamma(1-a_j+\alpha_js)}{\prod_{j=1}^{q}\Gamma(1-b_j+\beta_js)}}_{h_0(s)}\times\underbrace{\pi^{m+n-p}\frac{\prod_{j=n+1}^{p}\sin[\pi(a_j-\alpha_js)]}
{\prod_{j=1}^{m}\sin[\pi(b_j-\beta_js)]}}_{h_1(s)}.
\end{multline*}
The contour $C$ is the right loop separating the poles $s=(b_j+\nu)/\beta_j$, $j=1,\ldots,m$, $\nu=0,1,\ldots$ (lying to the right of $C$, hence inside the contour) from the poles $s=(a_j-1-\nu)/\alpha_j$, $j=1,\ldots,n$, $\nu=0,1,\ldots$ (lying to the left of $C$, hence outside the contour).  Denote
\begin{equation}\label{eq:mu-beta-defined}
\mu=\sum_{j=1}^{q}\beta_j-\sum_{j=1}^{p}\alpha_j, ~~~~\beta=\prod_{j=1}^{p}\alpha_{j}^{\alpha_{j}}\prod_{j=1}^{q}\beta_{j}^{-\beta_{j}},
\end{equation}
and
\begin{equation}\label{eq:delta0}
\delta_0=\left(\sum_{j=1}^{m}\beta_j-\sum_{j=n+1}^{p}\alpha_j\right)\pi.
\end{equation}
Next, by \cite[(4.26),(4.27),p.270]{Braaksma}
\begin{equation}\label{eq:Q-defined}
Q(z)=\sum\text{residues of }h(s)z^s\text{ in the points }s=(a_j-1-\nu)/\alpha_j, j=1,\ldots,n, \nu=0,1,\ldots
\end{equation}
\begin{equation}\label{eq:P-defined}
P(z)=\sum\text{residues of }h_0(s)z^s\text{ in the points }s=(a_j-1-\nu)/\alpha_j, j=1,\ldots,p, \nu=0,1,\ldots
\end{equation}
In \cite[Lemma~4, p.263]{Braaksma} Braaksma uses $\sin(z)=e^{iz}(1-e^{-2iz})/(2i)$ to expand:
\begin{multline*}
h_1(s)=(2\pi{i})^{m+n-p}\frac{\prod_{j=n+1}^{p}e^{i\pi(a_j-\alpha_js)}\left(1-e^{-2i\pi(a_j-\alpha_js)}\right)}
{\prod_{j=1}^{m}e^{i\pi(b_j-\beta_js)}\left(1-e^{-2i\pi(b_j-\beta_js)}\right)}
\\
=\underbrace{(2\pi{i})^{m+n-p}e^{i\pi\left\{\sum_{j=n+1}^{p}{a_j}-\sum_{j=1}^{m}b_j\right\}}}_{c_0}
\underbrace{e^{i\pi{s}\left\{\sum_{j=1}^{m}\beta_j-\sum_{j=n+1}^{p}{\alpha_j}\right\}}}_{=e^{i\gamma_0{s}}}
\frac{\prod_{j=n+1}^{p}\left(1-e^{2i\pi(\alpha_js-a_j)}\right)}
{\prod_{j=1}^{m}\left(1-e^{2i\pi(\beta_js-b_j)}\right)}.
\end{multline*}
For sufficiently large positive $\Im(s)$ we will have $\Re(2i\pi(\beta_js-b_j))<0$ so that $|e^{2i\pi(\beta_js-b_j)}|<1$, and we can expand
$\left(1-e^{2i\pi(\beta_js-b_j)}\right)^{-1}$ in geometric series:
$$
\frac{1}{1-e^{2i\pi(\beta_js-b_j)}}=\sum\limits_{\nu=0}^{\infty}e^{2i\pi\nu(\beta_js-b_j)}=\sum\limits_{\nu=0}^{\infty}\underbrace{e^{-2i\pi{\nu}b_j}}_{B_{j,\nu}}e^{2i\pi\nu\beta_js}.
$$
Expanding we get:
\begin{multline*}
\prod_{j=n+1}^{p}\left(1-e^{2i\pi(\alpha_js-a_j)}\right)=1-A_{n+1}e^{2i\pi\alpha_{n+1}s}-A_{n+2}e^{2i\pi\alpha_{n+2}s}-\ldots-A_{p}e^{2i\pi\alpha_{p}s}
\\
+\ldots+A_{n+1}A_{n+2}\cdots A_{p}e^{2i\pi(\alpha_{n+1}+\alpha_{n+2}+\cdots+\alpha_{p})s},
\end{multline*}
where $A_{k}=e^{-2i\pi{a_k}}$.  Multiplying all these expressions yields:
$$
h_1(s)=c_0e^{i\gamma_0{s}}\left(1+\hat{c}_1e^{i\gamma_1{s}}+\hat{c}_2e^{i\gamma_2{s}}+\cdots\right)
$$
with $0<\gamma_1<\gamma_2<\cdots$, where
\begin{equation}\label{eq:gamma0-defined}
\gamma_0=\pi\left\{\sum\nolimits_{j=1}^{m}\beta_j-\sum\nolimits_{j=n+1}^{p}{\alpha_j}\right\}
\end{equation}
and
\begin{equation}\label{eq:gamma1-defined}
\gamma_1=2\pi\min(\alpha_{n+1},\alpha_{n+2},\ldots,\alpha_{p},\beta_1,\ldots,\beta_m).
\end{equation}
If $\Im(s)$ is negative and of large absolute value, we can use $\sin(z)=-e^{-iz}(1-e^{2iz})/(2i)$ to similarly expand
$$
h_1(s)=d_0e^{-i\gamma_0{s}}\left(1+\hat{d}_1e^{-i\gamma_1{s}}+\hat{d}_2e^{-i\gamma_2{s}}+\cdots\right),
$$
where
$$
d_0=(-2\pi{i})^{m+n-p}e^{i\pi\left\{\sum\nolimits_{j=1}^{m}b_j-\sum\nolimits_{j=n+1}^{p}{a_j}\right\}}
=(2\pi)^{2(m+n-p)}/c_0,
$$
and $d_j=1/c_j$.  In \cite[Definition~I,p.265]{Braaksma} Braaksma takes the (unordered) multi-set
$$
\{\ldots,-\gamma_{2},-\gamma_{1},-\gamma_{0},\gamma_{0},\gamma_{1},\gamma_{2},\ldots\},
$$
sorts it in ascending order, removes the repeated elements and calls the new set
$$
\{\ldots,\delta_{-2},\delta_{-1},\delta_{0},\delta_{1},\delta_{2},\ldots\}~\text{with}~\delta_0=\gamma_0.
$$
\textbf{We next assume $\gamma_0=\delta_0=0$.}  Then $\delta_j=-\delta_{-j}$ and $\delta_j>0$ for $j=1,2,\ldots$, and, according to \cite[Definition~II,p.265]{Braaksma}, we see that $\kappa=0$ (by definition $\kappa$ is the index such that $\delta_{\kappa}=-\delta_0$).  Further, if an integer $r>0$ we will have case a) in \cite[Definition~II,p.265]{Braaksma} and $C_r=c_r$, $D_r=0$; if $r<0$ we have case b) and $C_r=0$, $D_r=-d_r$; if $r=0$, then we have case c) and $C_0=c_0$, $D_0=-d_0=-(2\pi)^{2(m+n-p)}/c_0$.  Note that the above formulas imply that (this is a slight modification of \cite[(4.10),p.266]{Braaksma})
\begin{equation}\label{eq:FourierUpper}
h_1(s)=\sum\limits_{j=0}^{\infty}c_je^{i\delta_{j}s},~~~c_j=c_0\hat{c}_j~\text{for}~j=1,2,\ldots
\end{equation}
for $\Im(s)>0$, and
\begin{equation}\label{eq:FourierLower}
h_1(s)=\sum\limits_{j=0}^{\infty}d_je^{-i\delta_{j}s}~~~c_j=d_0\hat{d}_j~\text{for}~j=1,2,\ldots
\end{equation}
for $\Im(s)<0$; each series obviously uniformly converges: first for $\Im(s)\ge\epsilon>0$, while the second for
$\Im(s)\le-\epsilon<0$. Hence, $c_j$, $d_j$ can be viewed as generalized Fourier coefficients in the upper and lower half-planes, while $\delta_j$ are ''generalized frequencies''.

Having made these preparations we can formulate (a part of) \cite[Theorem~3,p.280]{Braaksma} as follows.
\begin{theorem}\label{th:BraaksmaT2}
Suppose $\mu=\gamma_0=0$ \emph{(}see \eqref{eq:mu-beta-defined} and \eqref{eq:gamma0-defined} for definitions\emph{)}.  Then $H(z)$ can be continued analytically into the sector
$$
\Delta_{-}=\{z: -\gamma_1<\arg(z)<0\},
$$
where $\gamma_1$ is defined in \eqref{eq:gamma1-defined}, by
$$
H(z)=\frac{1}{2\pi{i}}\int\nolimits_{D}h_0(s)\left\{h_1(s)-c_0\right\}z^sds-V(z);
$$
and into the sector
$$
\Delta_{+}=\{z: 0<\arg(z)<\gamma_1\}
$$
by
$$
H(z)=\frac{1}{2\pi{i}}\int\nolimits_{D}h_0(s)\left\{h_1(s)-d_0\right\}z^sds-V(z);
$$
the contour $D$ starts at $s=-i\infty+\sigma$ and terminates at $s=i\infty+\sigma$, $\sigma$ is arbitrary real number,  leaves the points $s=(a_j-\nu-1)/\alpha_j$, $j=1,\ldots,p$, $\nu=0,1,\ldots$ \emph{(}the poles of $h_0(s)$\emph{)}  on the left and the points $s=(b_j+\nu)/\beta_j$, $j=1,\ldots,m$, $\nu=0,1,\ldots$, which do not coincide with  the points $(a_j-\nu-1)/\alpha_j$, $j=1,\ldots,p$, $\nu=0,1,\ldots$, \emph{(}the poles of $h(s)$ which are not poles of $h_0(s)$\emph{)} on the right.  The function $V(z)$ is the sum of residues of $h(s)z^s$ in the points which are poles of both $h(s)$ and $h_0(s)$.  This function is  equal to zero if
$$
\frac{b_j+\nu}{\beta_j}\ne\frac{a_h-1-\lambda}{\alpha_h}
$$
for all $j=1,\ldots,m$, $h=n+1,\ldots,p$, $\nu,\lambda=0,1,2,\ldots$

Moreover, the function $H(z)$ can be continued analytically from the sector $\Delta_{-}$ to the domain $|z|>\beta^{-1}$ by the formula
$$
H(z)=Q(z)-c_0P(z),
$$
and from the sector $\Delta_{+}$ to the domain $|z|>\beta^{-1}$ by the formula
$$
H(z)=Q(z)-d_0P(z).
$$
Here $Q$ and $P$ are defined in \eqref{eq:Q-defined} and \eqref{eq:P-defined}, respectively. Both series \eqref{eq:Q-defined} and \eqref{eq:P-defined} converge for
$|z|>\beta^{-1}$.
\end{theorem}

Our next goal is to find conditions that would guarantee that the analytic continuations to the sectors $\Delta_{-}$ and $\Delta_{+}$ are analytic in the ''large'' sector $\Delta_{-}\cup[0,\infty)\cup\Delta_{+}$. First, we have to require $c_0=d_0$.  Writing the definitions of $c_0$ and $d_0$ and performing some simple calculations we arrive at the condition:
$$
\sum_{j=n+1}^{p}{a_j}-\sum_{j=1}^{m}b_j=\frac{p-m-n}{2}+\eta,
$$
where $\eta$ is any integer. In this case we have
$$
c_0=d_0=(2\pi)^{m+n-p}e^{i\pi(m+n-p)/2}e^{i\pi(p-m-n)/2}e^{i\pi{\eta}}=(-1)^{\eta}(2\pi)^{m+n-p}.
$$
In particular, if $n=0$, $m=q$ and $\eta=1$ we obtain
$$
\sum_{j=1}^{p}{a_j}-\sum_{j=1}^{q}b_j=\frac{p-q}{2}+1,
$$
which corresponds to the case considered in \cite{KPCMFT2017}.  Now in this particular situation ($\mu=\delta_0=0$, $c_0=d_0$) we obtain a strengthened version of \cite[Lemma~4a, p.265]{Braaksma}:

\begin{lemma}\label{lm:4a}
Suppose $\mu=\delta_0=0$, $c_0=d_0$.  Then there exists a positive constant $K=K(\epsilon)$, independent of $s$, such that for any $\epsilon>0$
$$
|(h_1(s)-c_0)e^{-i\delta_{1}s}|\leq K~\text{if}~\Im(s)\ge\epsilon
$$
and
$$
|(h_1(s)-c_0)e^{i\delta_{1}s}|\leq K~\text{if}~\Im(s)\le-\epsilon.
$$
\end{lemma}
For the proof note that $c_0=d_0$ and use \eqref{eq:FourierUpper} and \eqref{eq:FourierLower} to get straightforward term-by-term estimates, similarly to the proof of \cite[Lemma~4a, p.265]{Braaksma}. Note, however, that we get a larger difference $2\delta_1$ between the exponents of the multipliers $e^{-i\delta_{1}s}$ and $e^{i\delta_{1}s}$ in the first and second bounds, than that in \cite[Lemma~4a, p.265]{Braaksma}.  This becomes possible by virtue of the condition $c_0=d_0$.

We now formulate an extension of a particular case of \cite[Theorem~2, p.280]{Braaksma} for the values of parameters $\mu=\delta_0=0$ and $c_0=d_0=(-1)^{\eta}(2\pi)^{m+n-p}$ ($\eta\in\Z$ is arbitrary).

\begin{theorem}\label{th:bigsector}
Suppose $\mu=\delta_0=0$, where $\mu$ and $\delta_0$ are defined in \eqref{eq:mu-beta-defined} and \eqref{eq:delta0}, respectively, and $c_0=d_0$. Then $c_0=d_0=(-1)^{\eta}(2\pi)^{m+n-p}$, where
$$
\eta=\sum_{j=n+1}^{p}{a_j}-\sum_{j=1}^{m}b_j-\frac{p-m-n}{2}\in\Z.
$$
Furthermore, the function $H^{m,n}_{p,q}(z)$ defined in \eqref{eq:Hmnpq} can be continued analytically to the sector
\begin{equation}\label{eq:bigsector}
\Delta_{\gamma_1}=\{z:-\gamma_1<\arg(z)<\gamma_1\},
\end{equation}
where $\gamma_1$ is  defined in \eqref{eq:gamma1-defined} by means of
\begin{equation}\label{eq:Hcontc0d0}
H(z)=\frac{1}{2\pi{i}}\int\nolimits_{D}h_0(s)\left\{h_1(s)-c_0\right\}z^sds-V(z);
\end{equation}
the contour $D$ and the function $V(z)$ have been defined in Theorem~\ref{th:BraaksmaT2}. Moreover, the function $H(z)$ can be continued analytically from this sector to the domain $|z|>\beta^{-1}$ on the Riemann surface of the logarithm  by the formula
$$
H(z)=Q(z)-c_0P(z),
$$
Here $Q$ and $P$ are defined in \eqref{eq:Q-defined} and \eqref{eq:P-defined}, respectively. Both series \eqref{eq:Q-defined} and \eqref{eq:P-defined} converge for $|z|>\beta^{-1}$.
\end{theorem}
\textbf{Proof.} We follow the proof of \cite[Theorem~2, p.280]{Braaksma} for $r=1$ and \cite[(6.10), p.280]{Braaksma} taking the form:
$$
f(s)=h_0(s)\beta^{-s}\{h_1(s)-c_0\}e^{-i\gamma_1s},
$$
so that, in view of the first bound in Lemma~\ref{lm:4a}, we can apply  \cite[Lemma~6, p.271]{Braaksma} with $\beta{z}e^{i\gamma_1}$ playing the role of $z$. Then \cite[Lemma~6, p.271]{Braaksma} applies for $\arg(z)>-\gamma_1$.  Instead of formula \cite[(6.11), p.280]{Braaksma} we take
$$
f(s)=h_0(s)\beta^{-s}\{h_1(s)-c_0\}e^{i\gamma_1s}.
$$
In view of the second bound in Lemma~\ref{lm:4a}, we can now apply  \cite[Lemma~6a, p.277]{Braaksma} with $\beta{z}e^{-i\gamma_1}$ playing the role of $z$ and hence \cite[Lemma~6a, p.277]{Braaksma} applies for $\arg(z)<\gamma_1$.
This leads to the validity of \cite[(6.12), p.280]{Braaksma} for $-\gamma_1<\arg(z)<\gamma_1$.  The remaining part of the proof of \cite[Theorem~2, p.280]{Braaksma} remains intact. $\hfill\square$

For the function $H^{q,0}_{p,q}(z)$ ($m=q$, $n=0$) we have $Q(z)=0$ by \eqref{eq:Q-defined}, so that if $\mu=0$ (and hence in this case $\delta_0=\pi\mu=0$) and
\begin{equation}\label{eq:eta}
\eta=\sum_{j=1}^{p}{a_j}-\sum_{j=1}^{q}b_j+\frac{q-p}{2}\in\Z,
\end{equation}
(note that $\eta=\alpha-1/2$, where $\alpha$ is defined in \cite[(3.24)]{Braaksma}), then $H^{q,0}_{p,q}(z)$ is analytic in the sector
$$
-\gamma_1<\arg(z)<\gamma_1,
$$
where $\gamma_1$ is  defined in \eqref{eq:gamma1-defined}, its analytic continuation to this sector is given by \eqref{eq:Hcontc0d0} (with $c_0=(-1)^{\eta}(2\pi)^{q-p}$) and from this sector to the domain  $|z|>\beta^{-1}$  by
$$
(-1)^{\eta+1}(2\pi)^{q-p}\sum\text{residues of }h_0(s)z^s\text{ in the points }s=(a_j-1-\nu)/\alpha_j, j\in\{1,\ldots,p\}, ~\nu\in\N_0,
$$
according to \eqref{eq:P-defined}.

\section{Expansion in the neighborhood of the singular point}
We now focus on the function $H^{q,0}_{p,q}(z)$  with $\mu=0$.  For further convenience we also apply the change of variable $s\to-s$ in the integral \eqref{eq:Hmnpq} to obtain for $|z|<\beta^{-1}$ the definition
\begin{equation}\label{eq:Hq0pq}
H^{q,0}_{p,q}\left(\!z\left|\begin{array}{l}(\alpha_1,a_1),\ldots(\alpha_p,a_p)\\(\beta_1,b_1),\ldots(\beta_q,b_q)\end{array}\right.\!\right)
=\frac{1}{2\pi{i}}\int\limits_{C'}\frac{\prod_{j=1}^{q}\Gamma(b_j+\beta_js)}{\prod_{j=1}^{p}\Gamma(a_j+\alpha_js)}z^{-s}ds,
\end{equation}
where the contour $C'=-C$, it is the left loop starting at $-\infty-ik$, terminating at $-\infty+ik$ and encompassing all the poles $s=(b_j-\nu)/\beta_j$, $j=1,\ldots,q$, $\nu=0,1,\ldots$ (so that $k>\Im(b_j/\beta_j)$, $j=1,\ldots,q$).  To be consistent with \cite{KPJMS2018,KPCMFT2017} denote
$$
\rho=\beta^{-1}>0.
$$
Assuming that $\eta=1$, we see by Theorem~\ref{th:bigsector} that the function $H^{q,0}_{p,q}(\rho{t})$  is analytic in the domain $G:=U_{-}\cup\Delta_{\gamma_1}$, where $U_{-}=\{t:|t|<1\}\setminus(-1,0]$ denotes the unit disk cut along the interval $(-1,0]$ (we need a cut in view of the branch point at $t=0$) and
$$
\Delta_{\gamma_1}=\{-\gamma_1<\arg(t)<\gamma_1\},~~~\gamma_1=2\pi\min(\alpha_1,\ldots,\alpha_p,\beta_1,\ldots,\beta_q),
$$
is the sector defined in \eqref{eq:bigsector} (note that change of variable $z\to\rho{t}$ does not alter this sector). Hence, for arbitrary  complex $\sigma$ we have the power series expansion
\begin{equation}\label{eq:Hqopqeta0}
\phi(t)=t^{-\sigma}H^{q,0}_{p,q}({\rho}t)=\sum\limits_{n=0}^{\infty}V_n(\sigma)(1-t)^n
\end{equation}
convergent in the disk $|1-t|<R$, where $R$ is the distance from $t=1$ to the boundary of the domain $G$.  Elementary geometry shows that $R=1$ if $\gamma_1\ge\pi/3$ (in view of \eqref{eq:gamma1-defined}, this is equivalent to $\min(\alpha_1,\ldots,\alpha_p,\beta_1,\ldots,\beta_q)\ge1/6$) and $R=2\sin(\gamma_1/2)$ if $\gamma_1<\pi/3$.
Assume for a moment that $\gamma_1\ge\pi/3$. Multiplying the above expansion by $t^{z-1}$  and integrating term by term  (to be justified below) we get
$$
\int_0^1 t^{z+\sigma-1}t^{-\sigma}H^{q,0}_{p,q}({\rho}t)dt=\sum\limits_{n=0}^{\infty}V_n(\sigma)\int_{0}^{1}t^{z+\sigma-1}(1-t)^ndt
=\sum\limits_{n=0}^{\infty}\frac{V_n(\sigma)n!}{(z+\sigma)_{n+1}}.
$$
On the other hand, by \cite[Theorem~6]{KPCMFT2016}
$$
\int_{0}^{1}t^{z-1}H^{q,0}_{p,q}({\rho}t)dt=\rho^{-z}\int_{0}^{\rho}u^{z-1}H^{q,0}_{p,q}(u)du
=\rho^{-z}\frac{\prod_{j=1}^{q}\Gamma(\beta_{j}z+b_j)}{\prod_{j=1}^{p}\Gamma(\alpha_{j}z+a_j)}
$$
or
\begin{equation}\label{eq:invfactorialgammas}
W(z)=\rho^{-z}\frac{\prod_{j=1}^{q}\Gamma(\beta_jz+b_j)}{\prod_{j=1}^{p}\Gamma(\alpha_jz+a_j)}
=\sum\limits_{n=0}^{\infty}\frac{V_n(\sigma)n!}{(z+\sigma)_{n+1}}.
\end{equation}
This inverse factorial series was derived in \cite[Theorem~1]{KPJMS2018}, where the coefficients were found explicitly.  However, the condition $\gamma_1\ge\pi/3$ was mistakenly omitted in this reference.

Here we give a more direct proof of this inverse factorial series expansion than that presented in \cite{KPJMS2018}, which will lead to a corrected proof of the expansion  \cite[Theorem~1]{KPCMFT2017} generalizing \eqref{eq:Hqopqeta0} to the case $\eta\ne1$.  Our proof will be based on a theorem due to N{\o}rlund.  To formulate this theorem  we need to define the Hadamard order of a function $f=\sum_{n\ge0}f_nx^n$ holomorphic inside the unit disk $|x|<1$.  According to \cite[23(1),p.46]{Norlund1926} the Hadamard order of $f$ on the circle $|x|=1$ is defined by:
\begin{equation}\label{eq:Hadamardorder}
\omega=\limsup\limits_{n\to\infty}\frac{\log|nf_n|}{\log(n)}.
\end{equation}

N{\o}rlund proved the following theorem (the formulation below is a combination of \cite[Theorem \S87,p.187]{Norlund1926} and \cite[Theorem \S88,p.189]{Norlund1926}).

\begin{theorem}\label{th:Norlund}
Suppose
$$
\phi(t)=\sum\limits_{s=0}^{\infty}a_s(1-t)^s
$$
is holomorphic for $|1-t|<1$ and has a finite order on the circle $|1-t|=1$.  Then
$$
F(w)=\int_0^1t^{w-1}\phi(t)dt
$$
can be expanded in the inverse factorial series
$$
F(w)=\sum\limits_{s=0}^{\infty}\frac{a_ss!}{(w)_{s+1}}
$$
convergent in a half-plane $\Re(w)>\lambda$ \emph{(}save the points $w=0,-1,-2,\ldots$\emph{)} and divergent for $\Re(w)<\lambda$. If the series $\sum_{s=0}^{\infty}a_s$ diverges, the convergence abscissa $\lambda=h-1$, where $h$ is the order of $t^{-1}\phi(t)$ on the the circle $|1-t|=1$; if $\sum_{s=0}^{\infty}a_s$ converges, then $\lambda=h'-1$, where $h'$ is the order of $[\phi(t)-\phi(0+)]/t$ on the circle $|1-t|=1$.
\end{theorem}

Assume that $\gamma_1>\pi/3$. Then, according to the above, for some $\varepsilon>0$ and any $0<\phi<\pi/2$ the function  $H^{q,0}_{p,q}(\rho{t})$ is holomorphic in the domain
$$
\Delta(\phi,\varepsilon)=\{t:|1-t|\le 1+\varepsilon~\text~{and}~-\pi+\phi\le\arg(t)\le\pi-\phi\},
$$
which is larger than the disk $|1-t|<1$.  We want to apply Theorem~\ref{th:Norlund} with $\phi(t)=t^{-\sigma}H^{q,0}_{p,q}(\rho{t})$ and $w=z+\sigma$.  Hence we need to determine the order of $\phi(t)t^{-1}=t^{-\sigma-1}H^{q,0}_{p,q}(\rho{t})$ on the circle $|1-t|=1$ and the order of $[\phi(t)-\phi(0+)]/t$ when $\phi(0+)\ne\infty$. First recall the shifting property \cite[(2.1.5)]{KilSaig}
$$
t^{-\sigma}H^{q,0}_{p,q}\left(\!\rho{t}\left|\begin{array}{l}(\alpha_1,a_1),\ldots(\alpha_p,a_p)\\(\beta_1,b_1),\ldots(\beta_q,b_q)\end{array}\right.\!\right)
=\rho^{\sigma}H^{q,0}_{p,q}\left(\!\rho{t}\left|\begin{array}{l}(\alpha_1,a_1-\sigma\alpha_1),\ldots(\alpha_p,a_p-\sigma\alpha_p)
\\(\beta_1,b_1-\sigma\beta_1),\ldots(\beta_q,b_q-\sigma\beta_q)\end{array}\right.\!\right)
$$
Next, suppose
$$
\mathcal{P}_{\sigma}=\left\{(s,r): s=\text{pole of}~h_{\sigma}(s); r=\text{its multiplicity}\right\},
$$
where $s$ runs over all distinct poles of the integrand of $t^{-\sigma}H^{q,0}_{p,q}(\rho{t})$ which is given by
$$
h_{\sigma}(s)=\frac{\prod_{j=1}^{q}\Gamma(b_j-\sigma\beta_j+\beta_js)}{\prod_{j=1}^{p}\Gamma(a_j-\sigma\alpha_j+\alpha_js)},
$$
and $r$ stands for the multiplicity of the corresponding pole.  We need the following definitions:  denote by  $\mathcal{P}'_{\sigma}$ the set $\mathcal{P}_{\sigma}$ where the elements of the form $(m,1)$, $m\in\Z_{\le0}$, have been removed (in other words simple poles at non-positive integers are removed); and define
\begin{equation}\label{eq:hatbeta}
\hat{\beta}(\sigma)=\max\{\Re(s): (s,r)\in\mathcal{P}'_{\sigma}\}.
\end{equation}
Note that $\hat{\beta}(\sigma)=\max\nolimits_{1\le{j}\le{q}}\Re(\sigma-b_j/\beta_j)$ except when this maximum is attained for $\sigma-b_j/\beta_j\in\Z_{\le0}$ and the corresponding pole is simpe. It is also convenient to write $\mathcal{P}_{\sigma,1}$ for the set of the first components  of the elements of $\mathcal{P}_{\sigma}$ (i.e. the set of poles of $h_{\sigma}(s)$ without multiplicities). Similarly for $\mathcal{P}'_{\sigma}$.

\begin{lemma}\label{lm:Hadamard}
Let $\phi(t)$ be defined in \eqref{eq:Hqopqeta0}.

If the series $\sum_{n\ge0}V_n(\sigma)$ diverges, then the Hadamard order $h$ of $t^{-1}\phi(t)$ on the circle $|1-t|=1$ is given by $h=\hat{\beta}(\sigma)+1$.

If the series $\sum_{n\ge0}V_n(\sigma)$ converges, then the Hadamard order $h'$ of $[\phi(t)-\phi(0+)]/t$
on the circle $|1-t|=1$  is given by the same expression $h'=\hat{\beta}(\sigma)+1$.
\end{lemma}
\textbf{Proof.} Denote
$$
I_{\sigma}:=-\mathcal{P}_{\sigma,1}\setminus\mathcal{P}'_{\sigma,1}.
$$
By definition $I_{\sigma}\subset\Z_{\ge0}$.  The residue theorem applied to the definition of $t^{-\sigma}H^{q,0}_{p,q}(\rho{t})$ gives as $t\to0$ \cite[Theorem~1.5]{KilSaig}:
\begin{multline}\label{eq:Hexpansion}
\phi(t)=t^{-\sigma}H^{q,0}_{p,q}(\rho{t})=\rho^{\sigma}\sum\mathop{\mathrm{res}}\limits_{s\in\mathcal{P}'_{\sigma,1}}h_{\sigma}(s)(\rho{t})^{-s}
+\rho^{\sigma}\sum\mathop{\mathrm{res}}\limits_{s\in\mathcal{P}_{\sigma,1}\setminus\mathcal{P}'_{\sigma,1}}h_{\sigma}(s)(\rho{t})^{-s}
\\
=(\log{t})^{r-1}\sum\limits_{k=1}^{l}C_kt^{-\hat{b}_{k}}\left(1+o(1)\right)+\sum\limits_{j\in I_{\sigma}\subset\Z_{\ge0}}A_jt^{j},
\end{multline}
where $\hat{b}_{k}=\sigma-b_{j_k}/\beta_{j_k}$, $k=1,2,\ldots,l$, are the elements of $\mathcal{P}'_{\sigma,1}$ with $\Re(\hat{b}_{k})=\hat{\beta}(\sigma)$ and having maximal multiplicity (which we denoted by $r$) among all poles with real part $\hat{\beta}(\sigma)$.  The precise values of the non-zero constants $C_1,\ldots,C_l$ can be found in \cite[Theorem~1.5]{KilSaig}, but they are immaterial for our purposes here.

To determine the order note first that if the series $\sum_{n\ge0}V_n(\sigma)$ converges, then
$\lim_{t\to0}\phi(t)<\infty$ by Abel's theorem. Hence, by \eqref{eq:Hexpansion} this series diverges if
$\hat{\beta}(\sigma)>0$ or $\hat{\beta}(\sigma)=0$ and $r>1$ or even when $\hat{\beta}(\sigma)=0$, $r=1$ and $\Im(\hat{b}_{k})\ne0$ for some $k\in\{1,\ldots,l\}$ (because no finite limit $\lim_{t\to0}\phi(t)$ exists in all these cases).  As $\hat{b}_{k}\in\mathcal{P}'_{\sigma,1}$ the case $\hat{b}_{k}=0$ is impossible by definition of $\mathcal{P}'_{\sigma}$. Hence, we conclude that $\sum_{n\ge0}V_n(\sigma)$ diverges
if $\hat{\beta}(\sigma)\ge0$.  To determine the of Hadamard order of $t^{-1}\phi(t)$ in these situations we apply \cite[Theorem~VI.4]{FSBook} with
$$
f(x)=(1-x)^{-\sigma-1}H^{q,0}_{p,q}(\rho(1-x))
$$
and $\zeta=1$.  This gives the following asymptotic relation for the power series coefficients $V_n(\sigma+1)$
\begin{equation}\label{eq:fnasymp}
V_n(\sigma+1)=A_0+\sum\limits_{k=1}^{l}D_kn^{\hat{b}_{k}}(\log{n})^{r-1}(1+o(1))+o(1)
\end{equation}
with some complex constants $D_k$ and $o(1)$ being different in each term.  By definition of the Hadamard order \eqref{eq:Hadamardorder} the above asymptotic formula leads immediately to $h=\hat{\beta}(\sigma)+1$.

Next, we assume that the series $\sum_{n\ge0}V_n(\sigma)$ converges, so that $\lim_{t\to0}\phi(t)<\infty$.
By the above argument this is only possible if $\hat{\beta}(\sigma)<0$.  Then, by \eqref{eq:Hexpansion} we have
$$
\lim_{t\to0}\phi(t)=A_0
$$
($A_0$ may vanish).  The last term in \eqref{eq:Hexpansion} for $t^{-1}(\phi(t)-A_0)$ takes the form
$$
\sum\limits_{j\in I_{\sigma}\setminus\{0\}\subset\Z_{>0}}A_jt^{j-1}
$$
which is holomorphic around $t=0$ and does not affect the order of  $t^{-1}(\phi(t)-A_0)$ (see \cite[section~26]{Norlund1926}).
Hence \eqref{eq:fnasymp} with $A_0$ removed holds true again, and the order $h'$ of $t^{-1}(\phi(t)-A_0)$ still equals $h'=\hat{\beta}(\sigma)+1$. $\hfill\square$

The coefficients $V_n(\sigma)$ have been computed in \cite[Corollary~1]{KPJMS2018} and \cite[Theorem~1]{KPCMFT2017}
via rearrangement of Poincar\'{e} asymptotics of $W(z)$ defined in \eqref{eq:invfactorialgammas} as $z\to\infty$ into inverse factorial series (this method dates back to Stirling). Putting these facts together we arrive at the following statements.

\begin{theorem}\label{th:Hexpansion1}
Suppose $\gamma_1=2\pi\min(\alpha_1,\ldots,\alpha_p,\beta_1,\ldots,\beta_q)>\pi/3$, $\mu=0$ and $\eta=1$, where $\mu$ and $\eta$ are defined in \eqref{eq:mu-beta-defined} and  \eqref{eq:eta}, respectively.  Then the function
$$
t^{-\sigma}H^{q,0}_{p,q}\left(\!\rho{t}\left|\begin{array}{l}(\alpha_1,a_1),\ldots(\alpha_p,a_p)\\(\beta_1,b_1),\ldots(\beta_q,b_q)\end{array}\right.\!\right)
$$
defined in \eqref{eq:Hq0pq}, analytic in the sector $-\gamma_1<\arg(z)<\gamma_1$ by Theorem~\ref{th:bigsector}, can be developed in convergent power series \eqref{eq:Hqopqeta0} with coefficients given by
\begin{equation}\label{eq:Vn}
V_n(\sigma)=(2\pi)^{(q-p)/2}\prod\nolimits_{k=1}^{q}\beta_k^{b_k-1/2}\prod\nolimits_{j=1}^{p}\alpha_k^{1/2-a_j}
\sum\limits_{k=0}^{n}\frac{(-1)^kl_{n-k}}{k!(n-k)!}\Be^{(n+1)}_{k}(1-\sigma),
\end{equation}
where $\rho=\beta^{-1}$ and  $\beta$ is defined in \eqref{eq:mu-beta-defined}, $\Be^{(n+1)}_{k}(\cdot)$ is Bernoulli-N{\o}rlund polynomial \cite[(1)]{Norlund61}.  The coefficients $l_r$ satisfy the recurrence relation \emph{(}with $l_0=1$\emph{)}
\begin{equation}\label{eq:lrqm}
l_r=\frac{1}{r}\sum\limits_{m=1}^r q_m l_{r-m},
~~~q_m=\frac{(-1)^{m+1}}{m+1}\left[\sum\limits_{k=1}^q\frac{\Be^{(1)}_{m+1}(b_k)}{\beta_k^m}
-\sum\limits_{j=1}^p\frac{\Be^{(1)}_{m+1}(a_j)}{\alpha_j^m}\right].
\end{equation}
\end{theorem}

The following theorem is a corrected and refined version of \cite[Theorem~5]{KPJMS2018}.

\begin{theorem}\label{th:inversefact1}
Under conditions of Theorem~\ref{th:Hexpansion1} the inverse factorial series \eqref{eq:invfactorialgammas} converges in the half plane $\Re(z)>\lambda$ \emph{(}save the points $w=0,-1,-2,\ldots$\emph{)}  and diverges in $\Re(z)<\lambda$.  The convergence abscissa $\lambda=\hat{\beta}(\sigma)-\Re(\sigma)$, where $\hat{\beta}(\sigma)$ is defined in \eqref{eq:hatbeta}.
\end{theorem}

\textbf{Remark}.  If  $\hat{\beta}(\sigma)=\max\nolimits_{1\le{j}\le{q}}\Re(\sigma-b_j/\beta_j)$ (which will be the case if this maximum is not attained for $\sigma-b_j/\beta_j\in\Z_{\le0}$), then, clearly,
$$
\lambda=-\min\limits_{1\le{j}\le{q}}\Re(b_j/\beta_j).
$$
In the boundary case $\gamma_1=\pi/3$ the function $\phi(t)$ has three singularities of the circle $|1-t|=1$ and our calculation of order in Lemma~\ref{lm:Hadamard} becomes inapplicable.  Hence, the convergence abscissa will, generally speaking, be different from that given in Theorem~\ref{th:inversefact1}.

To get rid of the restriction $\eta=1$ we applied the following simple trick  in \cite[Corollary~1]{KPJMS2018}.
Given $W(z)$ from \eqref{eq:invfactorialgammas} with $\eta\ne1$, consider
$$
W_1(z)=W(z)\frac{\Gamma(z+\theta+\eta)}{\Gamma(z+\theta+1)}
=\rho^{-z}\frac{\Gamma(z+\theta+\eta)\prod_{j=1}^{q}\Gamma(\beta_jz+b_j)}{\Gamma(z+\theta+1)\prod_{j=1}^{p}\Gamma(\alpha_jz+a_j)},
$$
where $\theta$ is arbitrary real number. Obviously, the parameter $\eta$ defined in (\ref{eq:eta}) computed for $W_1(z)$ takes the value unity: $\eta_1=1$.  Hence, we can apply \eqref{eq:invfactorialgammas} to $W_1(z)$. Taking $\sigma=\theta+\eta$ we obtain:
\begin{equation}\label{eq:Winvfactorial}
W(z)=W_1(z)\frac{\Gamma(z+\theta+1)}{\Gamma(z+\theta+\eta)}
=\sum\limits_{n=0}^{\infty}\frac{h_n\Gamma(z+\theta+1)}{(z+\theta+\eta)_{n+1}\Gamma(z+\theta+\eta)}
=\sum\limits_{n=0}^{\infty}\frac{h_n\Gamma(z+\theta+1)}{\Gamma(z+\theta+\eta+n+1)},
\end{equation}
where the coefficients $h_n=V_n'(\theta+\eta)n!$ are obtained by applying formulas \eqref{eq:Vn}, \eqref{eq:lrqm} to the function $W_1(z)$. In other words, in these formulas we have to make the following replacements:
\begin{equation}\label{eq:substitute}
\begin{split}
&(\alpha_1,\ldots,\alpha_p)~\to~(\alpha_1,\ldots,\alpha_p,1),
~~~~~(a_1,\ldots,a_p)~\to~(a_1,\ldots,a_p,\theta+1),
\\
&(\beta_1,\ldots,\beta_q)~\to~(\beta_1,\ldots,\beta_q,1),
~~~~~(b_1,\ldots,b_q)~\to~(b_1,\ldots,b_q,\theta+\eta).
\end{split}
\end{equation}
The resulting expressions are given explicitly in \cite[Corollary~1]{KPJMS2018}.  The convergence abscissa of the above expansion equals that for $W_1(z)$ and hence by Theorem~\ref{th:inversefact1} is computed as follows:  remove non-positive integers from the set of simple  poles of $W_1(z)$ and denote by $\hat{\beta}_1(\theta+\eta)$ the maximum of the real parts of the remaining poles.  The convergence abscissa then equals  $\hat{\beta}_1(\theta+\eta)-\theta+\Re(\eta)$.

Given the convergent inverse factorial series (\ref{eq:Winvfactorial}) we can substitute the beta integral
$$
\frac{\Gamma(z+\theta+1)\Gamma(\eta+n)}{\Gamma(z+\theta+\eta+n+1)}
=\int_0^1t^{z+\theta}(1-t)^{n+\eta-1}
$$
and exchange the order of summation in integration. Then uniqueness of the inverse truncated Mellin transform immediately leads to
\begin{equation}\label{eq:Hexp-any-eta}
H^{q,0}_{p,q}\left(\!\rho{t}\left|\begin{array}{l}(\alpha_1,a_1),\ldots,(\alpha_p,a_p)\\(\beta_1,b_1),\ldots,(\beta_q,b_q)\end{array}\right.\!\right)
=t^{\theta+1}(1-t)^{\eta-1}\sum_{n=0}^{\infty}\frac{h_n}{\Gamma(\eta+n)}(1-t)^{n}.
\end{equation}
These steps have been rigourously justified in \cite[Theorem~1]{KPCMFT2017}, where we also presented alternative formulas for the coefficients.  Finally we remark that the radius of convergence of \eqref{eq:Hexp-any-eta} is the same as for the expansion \eqref{eq:Hqopqeta0}, where $\eta=1$. Indeed, this series has the same radius of convergence as the series $\sum\nolimits_{n\ge0}h_n(1-t)^{n}/n!$ which is an instance of \eqref{eq:Hqopqeta0} with substitutions \eqref{eq:substitute}.
Altogether, this yields the following statement generalizing Theorem~\ref{th:Hexpansion1} to arbitrary values of parameters.

\begin{theorem}\label{th:Hgeneral}
For arbitrary positive $\alpha_i$, $\beta_j$ and complex $a_i$, $b_j$, the series in \eqref{eq:Hexp-any-eta} with coefficients $h_n=V_n'(\theta+\eta)n!$, where $V_n'(\theta+\eta)$ are  given by \eqref{eq:Vn}, \eqref{eq:lrqm} with substitutions \eqref{eq:substitute}, has the radius of convergence $R=1$ if
$$
\gamma_1=2\pi\min(\alpha_1,\ldots,\alpha_p,\beta_1,\ldots,\beta_q)\ge\pi/3
$$
and $R=2\sin(\gamma_1/2)$ if $\gamma_1<\pi/3$.
\end{theorem}

In case $\gamma_1\le1/3$ we can easily get an expansion convergent in the disk of radius $1$, but at the price of having powers of $1-t^{\omega}$ instead of powers of $1-t$ with some $\omega>0$ .  Indeed, the substitution $t=u^{1/\omega}$ with $\omega>1$ extends the sector $\Delta_{\gamma_1}$ to the sector $\Delta_{\omega\gamma_1}$, where $\omega\gamma_1$ can be made greater than $\pi/3$ by choosing the appropriate $\omega$. Hence, $u\to H^{q,0}_{p,q}(\rho{u^{1/\omega}})$ is analytic in $|1-u|<1$ for such $\omega$.  On the other hand, according to \cite[(2.1.4)]{KilSaig}
$$
H^{q,0}_{p,q}\left(\!\rho{u^{1/\omega}}\left|\begin{array}{l}(\alpha_1,a_1),\ldots,(\alpha_p,a_p)\\(\beta_1,b_1),\ldots,(\beta_q,b_q)\end{array}\right.\!\right)
={\omega}H^{q,0}_{p,q}\left(\!\rho^{\omega}u\left|\begin{array}{l}(\omega\alpha_1,a_1),\ldots,(\omega\alpha_p,a_p)\\(\omega\beta_1,b_1),\ldots,(\omega\beta_q,b_q)\end{array}\right.\!\right).
$$
Writing $\hat{\alpha}_i=\omega\alpha_i$ and $\hat{\beta}_j=\omega\beta_j$, formula \eqref{eq:mu-beta-defined} and the definition $\rho=\beta^{-1}$ imply that $\hat{\rho}=\rho^{\omega}$, in view of  $\sum\nolimits_{i=1}^{p}\alpha_i=\sum\nolimits_{j=1}^{q}\beta_j$.  Hence, we can apply Theorem~\ref{th:Hgeneral} to the function on the right hand side of the above relation.  Returning to the original variable $t=u^{1/\omega}$, an extended version of  expansion \eqref{eq:Hexp-any-eta} then takes the form:
\begin{equation}\label{eq:Hexp-any-omega}
H^{q,0}_{p,q}\left(\!\rho{t}\left|\begin{array}{l}(\alpha_1,a_1),\ldots,(\alpha_p,a_p)\\(\beta_1,b_1),\ldots,(\beta_q,b_q)\end{array}\right.\!\right)
={\omega}t^{\omega(\theta+1)}(1-t^{\omega})^{\hat{\eta}-1}
\sum_{n=0}^{\infty}\frac{\hat{h}_n}{\Gamma(\hat{\eta}+n)}(1-t^{\omega})^{n},
\end{equation}
where, obviously, the coefficients are calculated as described in  Theorem~\ref{th:Hgeneral}, but writing $\hat{\alpha}_i$, $\hat{\beta}_j$ instead of $\alpha_i$, $\beta_j$.  The series \eqref{eq:Hexp-any-omega} converges in the domain $|1-t^{\omega}|<1$.

\end{document}